\documentclass[12pt]{amsart}
\usepackage{graphicx}
\usepackage{color}
\usepackage{amsmath,amssymb,amsthm}
\usepackage[top=1in, bottom=1in, left=1in, right=1in]{geometry} 
\newcommand{\inprod}[2]{\langle #1,#2 \rangle}

\newcommand{\fpart}[2]{\tfrac{\partial #1}{\partial #2}}
\newcommand{\spart}[2]{\tfrac{\partial^2 #1}{\partial #2^2}}

\renewcommand{\(}{\left(}
\renewcommand{\)}{\right)}
\renewcommand{\epsilon}{\varepsilon}

\makeatletter
\newtheorem*{rep@theorem}{\rep@title}
\newcommand{\newreptheorem}[2]{%
\newenvironment{rep#1}[1]{%
 \def\rep@title{#2 \ref{##1}}%
 \begin{rep@theorem}}%
 {\end{rep@theorem}}}
\makeatother

\newtheorem{theorem}{Theorem}
\newtheorem{lemma}[theorem]{Lemma}

\newreptheorem{theorem}{Theorem}
\newtheorem{example}{Example}

\begin{document}


\title{Scaling limits of a model for selection at two scales}
\author{Shishi Luo, Jonathan C. Mattingly}
\date{\today}
\maketitle
\vspace{1cm}

\section*{Abstract}
The dynamics of a population undergoing selection is a central
topic in evolutionary biology. This question is particularly
intriguing in the case where selective forces act in opposing
directions at two population scales. For example, a fast-replicating
virus strain out-competes slower-replicating strains at the within-host scale. 
However, if the fast-replicating strain causes host morbidity and is
less frequently transmitted, it can be outcompeted by slower-replicating
 strains at the between-host scale. Here we consider a stochastic
ball-and-urn process which models this type of phenomenon.  We prove
the weak convergence of this process under two natural
scalings. The first scaling leads to a deterministic nonlinear
integro-partial differential equation on the interval $[0,1]$ with
dependence on a single parameter, $\lambda$. We show that the fixed points
of this differential equation are Beta distributions and that their
stability depends on $\lambda$ and the behavior of the initial data
around $1$. The second scaling leads to a measure-valued Fleming-Viot
process, an infinite dimensional stochastic process that is frequently associated with a
population genetics.

\section{Introduction}
We study the model, introduced in \cite{Luo:2014}, of a trait that is
advantageous at a local or individual level but disadvantageous at a
larger scale or group level. For example, an infectious virus strain
that replicates rapidly within its host will outcompete other virus
strains in the host. However, if infection with a heavy viral load is
incapacitating and prevents the host from transmitting the virus, the
rapidly replicating strain may not be as prevalent in the overall host
population as a slow replicating strain.

\begin{figure}
\includegraphics[width=.9\textwidth]{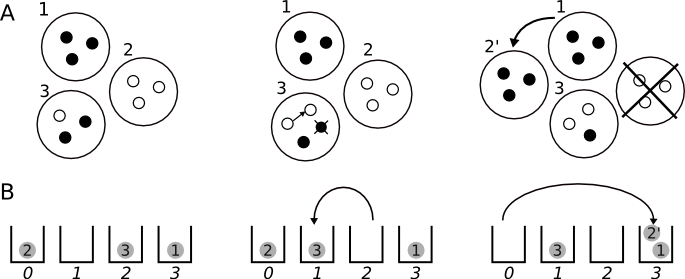}
\centering
\caption{Schematic of the particle process. (a) \emph{Left:} A
  population of $m=3$ groups, each with $n=3$ individuals of either
  type G (filled small circles) or type I (open small
  circles). \emph{Middle:} A type I individual replicates in group 3
  and a type G individual is chosen uniformly at random from group 3
  to die. \emph{Right: } Group 1 replicates and produces group
  $2'$. Group 2 is chosen uniformly at random to die. (b) The states
  in (a) mapped to a particle process. \emph{Left}: Group 2 has no
  type G individuals, represented by ball 2 in urn 0. Similarly, group
  3 is represented by ball 3 in urn 2 and group 1 by ball 1 in urn
  3. \emph{Middle}: The number of type G individuals in group 3
  decreases from two to one, therefore ball 3 moves to urn
  1. \emph{Right}: A group with zero type G individuals dies, while a
  group with three type G individuals is born. Therefore ball 2 leaves
  urn 0 and appears in urn 3 as ball $2'$.}
\label{fig_schematic}
\end{figure}

A simple mathematical formulation of this phenomenon is as
follows. Consider a population of $m\in\mathbb{N}$ groups. Each group
contains $n\in\mathbb{N}$ individuals. There are two types of
individuals: type I individuals are selectively advantageous at the
individual (I) level and type G individuals are selectively
advantageous at the group (G) level.  Replication and selection occur
concurrently at the individual and group level according to the Moran
process \cite{Durrett:2008} and are illustrated in
Fig~\ref{fig_schematic}. Type I individuals replicate at rate $1+s$,
$s\geq 0$ and type G individuals at rate $1$. When an individual gives
birth, another individual in the same group is selected uniformly at
random to die. To reflect the antagonism at the higher level of
selection, groups replicate at a rate which increases with the number
of type G indivduals they contain. As a simple case, we take this rate
to be $w(1+r \frac{k}{n})$, where $\frac{k}{n}$ is the fraction of
indivduals in the group that are type G, $r\geq 0$ is the selection
coefficient at the group level, and $w>0$ is the ratio of the rate of
group-level events to the rate of individual-level events. More
general functions for the group replication rate are possible, though
the subsequent analysis of the model may be less tractable. As with
the individual level, the population of groups is maintained at $m$ by
selecting a group uniformly at random to die whenever a group
replicates.  The offspring of groups are assumed to be identical to
their parent.

As illustrated in Fig~1, this two-level process is equivalent to a
ball-and-urn or particle process, where each particle represents a
group and its position corresponds to the number of type G individuals
that are in it.

Let $X_t^i$ be the number of type G individuals in group $i$ at time $t$. Then
\[\mu^{m,n}_t:=\frac{1}{m}\sum_{i=1}^m \delta_{X^i_t/n}\]
is the empirical measure at time $t$ for a given number of groups
$m$ and individuals per group $n$. $\delta_x(y) = 1$ if $x=y$ and zero
otherwise. The $X^i_t$ are divided by $n$ so that $\mu^{m,n}_t$ is a
probability measure on $E_n:=\{0,\frac{1}{n},\ldots,1\}$.

For fixed $T>0$, $\mu^{m,n}_t\in D([0,T],\mathcal{P}(E_n))$, the set
of c\`adl\`ag processes on $[0,T]$ taking values in
$\mathcal{P}(E_n)$, where $\mathcal{P}(S)$ is the set of probability
measures on a set $S$. With the particle process described above,
$\mu^{m,n}_t$ has generator
\begin{align}
  \label{gen_eqn}({L}^{m,n}\psi)(v) =
  \sum_{i,j}({R_1+wR_2})(v,v_{ij})[\psi(v_{ij})-\psi(v)]
\end{align}
where $v_{ij}:=v+\tfrac{1}{m}(\delta_{\frac{j}{n}}-\delta_{\frac{i}{n}})$, $\psi\in C_b(\mathcal{P}([0,1]))$ are bounded continuous functions,
and $v\in\mathcal{P}(E_n)\subset \mathcal{P}([0,1])$.  The transition rates $({R_1+wR_2})$ are given by
\[{R_1}\left(v,v_{ij}\right) = \left\{\begin{array}{lr} mv(\tfrac{i}{n})i\left(1-\tfrac{i}{n}\right)(1+s)& \text{ if }j=i-1, i<n\\
    mv(\tfrac{i}{n})i\left(1-\tfrac{i}{n}\right) & \text{ if } j=i+1, i>0\\
    0& \text{ otherwise }\end{array}\right.\]
    and
\[{R_2}\left(v,v_{ij}\right) = mv(\tfrac{i}{n})v(\tfrac{j}{n})(1+r\tfrac{j}{n}).\]
$R_1$ represents individual-level events while $R_2$ represents group-level events.
\\\\
\noindent\textbf{Acknowledgements:} 
The authors would like to thank Mike Reed and Katia Koelle for their 
roles in the collaboration out of which this paper's central model grew. We
would also like to thank Rick Durrett of a number of useful
discussions. JCM would like to thank the NSF for its support though DMS-08-54879.
SL would like to thank support from the NSF (grants NSF-EF-08-27416 and
DMS-0942760), NIH (grant R01-GM094402), and the Simons Institute for
the Theory of Computing.

\section{Main results}
We prove the weak convergence of this measure-valued process as
$m,n\to\infty$ under two natural scalings. The first scaling leads to
a deterministic partial differential equation. We derive a closed-form
expression for the solution of this equation and study its
steady-state behavior. The second scaling leads to an infinite
dimensional stochastic process, namely a Fleming-Viot process.

Let us briefly introduce some notation. By $m,n\to\infty$ we mean a
sequence $\{(m_k,n_k)\}_k$ such that for any $N$, there is an $n_0$
such that if $k\geq n_0$, $m_k, n_k\geq N$. We define $\inprod{f}{v} =
\int_0^1 f(x) v(dx)$ where $f$ is a test function and $v$ a
measure. Lastly, $\delta_x$ will denote the delta measure for both
continuous and discrete state spaces.

To provide intuition for the two scalings and the corresponding
limits, take $\psi$ to be of the form $\psi(v) = F(\inprod{g}{v})$,
where $g$ is some suitable function on $[0,1]$, and apply the
generator in \eqref{gen_eqn} to it:
\begin{align*}
L^{m,n}\psi(\mu) & = F'\cdot\Big\{\sum\left[\tfrac{1}{n}g''(\tfrac{i}{n})-sg'(\tfrac{i}{n})\right]\tfrac{i}{n}(1-\tfrac{i}{n})\mu(\tfrac{i}{n})\\
& \qquad\qquad +wr\left[\sum\tfrac{i}{n}g(\tfrac{i}{n})\mu(\tfrac{i}{n})-\sum g(\tfrac{i}{n})\mu(\tfrac{i}{n})\sum\tfrac{j}{n}\mu(\tfrac{j}{n})\right]\Big\}\\
& \qquad + \tfrac{1}{m}wF''\cdot\Big\{\sum g(\tfrac{i}{n})^2\mu(\tfrac{i}{n})-\left(\sum g(\tfrac{i}{n})\mu(\tfrac{i}{n})\right)^2\\
&\qquad\qquad\qquad\qquad + \tfrac{1}{2} r \sum(g(\tfrac{j}{n})-g(\tfrac{i}{n})^2\tfrac{j}{n}\mu(\tfrac{j}{n})\mu(\tfrac{i}{n})\Big\}\\
&\qquad +o(\tfrac{1}{m})+o(\tfrac{1}{n})
\end{align*}

This suggests two natural scalings. The first is to take
$m,n\to\infty$ without rescaling any parameters. The $g''$ and $F''$
terms vanish and we have a deterministic process. 
%
The precise statement of the weak convergence of the
finite state space system to the deterministic limit is in terms of a weak measure-valued solution to a partial differential equation:
\begin{theorem}
\label{det_thm}
Suppose the particles in the system described by $\mu^{m,n}_t$ are
initially independently and identically distributed according to the
measure $\mu^{m,n}_0$, where $\mu^{m,n}_0\to
\mu_0\in\mathcal{P}([0,1])$ as $m,n\to\infty$. Then, as $m,
n\to\infty$, $\mu^{m,n}_t\to\mu_t\in D([0,T], \mathcal{P}([0,1]))$
weakly, where $\mu_t$ solves the differential equation
\begin{align}
  \label{det_pde} \tfrac{d}{dt}\inprod{f}{\mu_t} = -\inprod{x(1-x)f'}{
    \mu_t}+\lambda\left[\inprod{xf}{\mu_t}-
    \inprod{f}{\mu_t}\inprod{x}{\mu_t}\right]
\end{align}
for any positive-valued test function $f\in C^1([0,1])$ and with
initial condition $\inprod{f}{\mu_0}$. Here, $\lambda:=\frac{wr}{s}$
and time has been sped up by a factor of $s$.
\end{theorem}

Throughout we will denote the measure-valued solutions to 
\eqref{det_pde} by $\mu_t(dx)$. We note that strong, density-valued solutions, denoted by 
$\eta_t(x)$, solve:
\begin{align}
\label{nondim_strong_pde}
\fpart{}{t}\eta_t= \fpart{}{x}\left[x(1-x)\eta_t\right] +
\lambda\eta_t\cdot\(x-\int_0^1y\eta_t(y)dy\)
\end{align} 
with initial density $\eta_0(x)$. In this more transparent form one can see
 that the first term on the right is a
flux term that transports density towards $x=0$ whereas the second
term is a forcing term that increases the density at values of $x$
above the mean of the density. The flux corresponds to the individual-level
moves: nearest neighbor moves in the particle system. The forcing term
corresponds to group-level moves: moves to occupied sites in the
particle system. 

We will see that if we start with an initial measure $\mu_0$ which is the sum
of delta measures, then the solution $\mu_t$ retains the same
form. More explicitly, if
\begin{align*}
  \mu_0(dx) = \sum_{i} a_i(0) \delta_{x_i(0)}(dx)
\end{align*}
where $x_i(0) \in [0,1]$, $a_i(0) >0$, and $\sum  a_i(0) =1$, then we
will see (from Lemma~\ref{method_char}) that the solution $\mu_t$ to
\eqref{det_pde} has the form
\begin{align*}
  \mu_t(dx) = \sum_{i} a_i(t) \delta_{x_i(t)}(dx)\,.
\end{align*}
Moreover, the parameters $(a_i(t),x_i(t))$ satisfy the following set
of coupled equations
\begin{equation}
 \left\{ \begin{aligned}
      \frac{d x_i}{dt}&= - x_i(1-x_i)\\
   \frac{d a_i}{dt}&= \lambda a_i \Big(  x_i -\inprod{y}{\mu_t} \Big)=\lambda a_i \Big(  x_i - \sum_j a_j x_j \Big)\,.
  \end{aligned}\right.\label{eq:particles}
\end{equation}
Notice that the positions of the delta masses change according to a
negative logistic function, independently of the other masses and the
density. The weight $a_i$ increases at time $t$ if the position of the particle
$x_i$ is above the mean, $\sum a_j x_j$, and decreases if it is
below the mean. To build intuition, it is instructive to consider some simple
examples of this form.

\begin{example}
According to \eqref{eq:particles}, if $\mu_0=\delta_1$, then
$\mu_t=\mu_0$. This can also be seen directly from \eqref{det_pde}.
  In the case of an initial condition containing some delta mass
at $1$, all of the rest of the mass will migrate towards
zero. Eventually all of the mass will be below the mean as the mass at
one will not move and will ever be increasing its mass as it is
always above the mean. Once this happens it is clear that all of the
mass will drain from all of the points not at one and hence $\mu_t
\rightarrow \delta_1$ as $t \rightarrow \infty$. This reasoning holds
in a more general setting and is included in Theorem~\ref{thm:steadyState}.
\end{example}

\begin{example}
According to \eqref{eq:particles}, if $\mu_0=\delta_0$, then
$\mu_t=\mu_0$. This too can be seen directly from \eqref{det_pde}.
In the case of an initial condition containing no mass at one and only 
finite number of masses total, the mass will eventually all move 
towards zero and hence  hence $\mu_t 
\rightarrow \delta_0$ as $t \rightarrow \infty$. If an infinite 
number of masses are allowed the situation is not as simple.
Theorem~\ref{thm:steadyState} hints at the possible complications by
giving an example of a density which is invariant.
\end{example}

Though $\delta_0$ is a fixed point of the system attracting many
initial configurations, it is not Lyapunov stable. This means that
even small perturbations of $\delta_0$ can lead to an arbitrary large
excursion away from $\delta_0$ even though the system eventually
returns to $\delta_0$. Rather than making a precise statement which
would require quantifying the size of a perturbation, consider the
example of $\mu_0 =(1-\epsilon)\delta_0 + \epsilon
\delta_{1-\alpha}$. As $\epsilon \rightarrow 0$, the distance between
$\mu_0$ and $\delta_0$ goes to zero in any reasonable metric. If we
write $\mu_t = (1-a_t)\delta_0 + a_t \delta_{x_t}$ then as $\alpha \rightarrow
0$ one can ensure that the system spends arbitrarily long time with
$x_t > \frac12$ and hence $a_t$ will grow to as close to one as one
wants in this time. Thus the system could be described as making an an
arbitrarily big excursion away from $\delta_0$  even though 
$\mu_t \rightarrow \delta_0$ as $t \rightarrow
\infty$.

It natural to ask if there are other fixed points beyond  $\delta_0$
and $\delta_1$. 
\begin{lemma}[Fixed points]\label{lem:fixedPoints}
 The measures delta $\delta_0$, $\delta_1$, and
  densities in the $Beta(\lambda-\alpha, \alpha)$ family of
  distributions:
\[\frac{1}{B(\lambda-\alpha, \alpha)}x^{\lambda-\alpha-1}(1-x)^{\alpha-1}\]
with $\alpha \in (0,\lambda)$, are fixed points of \eqref{det_pde}.
$B(\lambda-\alpha, \alpha)$ is the normalizing constant that makes the
density integrate to 1 over the interval $[0,1]$.
\end{lemma}

For measure-valued initial data, we show that the
basins of attraction for the fixed points are determined by whether
they charge the point $x=1$ and their
H\"older exponent around $x=1$.

\begin{theorem}[Steady state behavior] \label{thm:SteadyState}Consider measure valued
  solution $\mu_t(dx)$ to \eqref{det_pde} with initial probability
  measure $\mu_0(dx)$. 
 If $\mu_0(\{1\}) >0$ then
  \begin{align*}
\mu_t \rightarrow \delta_1\qquad\text{as}\qquad t 
  \rightarrow \infty 
  \end{align*}
and if $\mu_0([1-\epsilon,1]) =0$ for some $\epsilon >0$ then
\begin{align*}
\mu_t \rightarrow \delta_0\qquad\text{as}\qquad t 
  \rightarrow \infty \,.
  \end{align*}
Alternatively, suppose that for some $\alpha >0$ and $C>0$
  \begin{align*}
   x^{-\alpha} \mu_0([1-x,1]) \rightarrow C\qquad\text{as}\qquad x \rightarrow 0\,.
  \end{align*}
  If $\alpha <\lambda$, then 
  \begin{align*}
    \mu_t(dx) \rightarrow \text{Beta}(\lambda-\alpha, \alpha)\qquad\text{as}\qquad t\rightarrow\infty\,.
  \end{align*}
Otherwise, if $\alpha \geq  \lambda$, 
\begin{align*}
  \mu_t(dx)\to\delta_0(dx)\qquad\text{as}\qquad t\to\infty
\end{align*}\label{thm:steadyState}
\end{theorem}
%

The results of Theorem~\ref{thm:SteadyState} should be contrasted with
the original Markov chain before taking the limit $m,n \rightarrow \infty$. In the
Markov chain, all individuals eventually become either entirely type
G or type I. These two homogeneous states are absorbing states for
the individual level dynamics. The population level state made of individuals that are
all either homogeneous of type G or I is absorbing for the group
level dynamics.
Hence, the state of the system eventually becomes
composed entirely of homogeneous groups of solely G or
I and stays in that state for all future times.

These two absorbing states of the Markov chain, with finite $m$ and $n$, correspond to the states $\delta_0$ and
$\delta_1$ in the scaling limit. Hence the natural discretization for the
Beta distribution to the lattice $\{\frac{k}n : 0<k< n\}$, given by
\begin{equation*}
  \frac1{Z(m,n,\lambda,\alpha)}
  \big(\tfrac{k}{n}\big)^{\lambda-\alpha-1}\big(1-\tfrac{k}{n}\big)^{\alpha-1} ,
\end{equation*}
cannot be
invariant. (Here $Z$ is the normalization  constant which ensures the
probabilities sum to one.) However for large $m$ and $n$, it is reasonable to expect it
to be nearly invariant in the sense that if  the initial states
$\{X_i(0) : 1\leq i \leq m\}$ are independent and distributed as  the
discrete Beta distribution then the Markov chain dynamics will keep the
distribution close to the product of discretized Beta distributions
for a long time. The expectation of this time will grow to infinity as
$m,n \rightarrow  \infty$.

We will not pursue a rigorous proof of this near or quasi
invariance here. Nonetheless, we now briefly sketch the argument as we
understand it, giving the central points. If the distribution of the
Markov chain is close to a product of discretized Beta distributions,
then the empirical mean will be highly concentrated around the mean of
continuous Beta when $m$ and $n$ are large. Hence the generator
projected on to any $X_i$ is nearly decoupled from the other particles
and close to being Markovian. More precisely, the dynamics of any
fixed $X_i$ is well approximated in this setting by the
one-dimensional Markov chain obtained by replacing the mean of the
empirical measure in the full generator with the mean of the Beta
distribution. It is straightforward to see that for $m$ and $n$ large the
discretized Beta distribution is an approximate left-eigenfunction of
this one-dimensional generator with an eigenvalue which goes to zero
as $m,n \rightarrow \infty$.

All of these observations can be combined to show  that if the systems
starts in the product discretized Beta distribution then it will say
close to the product discretized Beta distribution for a long time if
$m$ and $n$ are large.

We now turn to the second scaling. Let $s=\frac{\sigma}{n}$,
$r=\frac{\rho}{m}$, and $\frac{n}{m}\to \theta$, and let $\nu^{m,n}_t$
denote the empirical measure under this scaling. The terms $F''$ and
$g''$ in the generator \eqref{gen_eqn} no longer vanish and the
process converges to a limit that is stochastic.  
Our weak convergence result is proved and stated in terms of a martingale problem.
\begin{theorem}
\label{fv_thm}
Suppose $\frac{n}{m}\to\theta$, $w = O(1)$, $s = \tfrac{\sigma}{n}$,
$r= \tfrac{\rho}{m}$, and we speed up time by a factor of $n$. Suppose
the particles in the rescaled $\nu^{m,n}_t$ process are initially
independently and identically distributed according to the measure
$\nu^{m,n}_0$ where $\nu^{m,n}_0\to \nu_0$ as $m,n\to\infty$. Then the
rescaled process converges weakly to $\nu_t$ as $m,n\to\infty$, where
$\nu_t$ satisfies the following martingale problem:
\begin{align}\label{fv_drift}
N_t(f) & = \inprod{f}{\nu_t} -\inprod{f}{\nu_0}- \int_0^t\inprod{Af}{\nu_{z}}dz\\
\notag&\qquad - w\theta \rho\int_0^t\left\{\int_0^1\int_0^1 f(x)V(z, \nu_{z}, y) Q(\nu_{z};dx,dy) \right\}dz
\end{align}
is a martingale with conditional quadratic variation
\begin{align}\label{fv_quad}
\langle N(f)\rangle_t&= 2w\theta \int_0^t\left[\int_0^1\int_0^1 f(x)f(y) Q(\nu_{\tau}; dx, dy)\right]d\tau
\end{align}
where
\begin{align*}
Af(x) &= x(1-x)\left[\tfrac{d^2}{dx^2}f(x) -\sigma \tfrac{d}{dx}f(x)\right]\\
V(t, \nu, x) &= x\\
Q(\nu; dx, dy) & = \nu(dx)(\delta_x(dy)-\nu(dy))
\end{align*}
and $f\in C^2([0,1])$.
\end{theorem}

The drift part of the martingale \eqref{fv_drift} comprises a second order partial
differential operator $A$ and the centering term from the global jump
dynamics (the expression in curly brackets). The entire process is a Fleming-Viot process \cite{Fleming:1979}. Fleming-Viot processes
frequently arise in models of population genetics (see
\cite{Ethier:1993} for a review). In these contexts, the variable $x$
can represent the geographical location of an individual, or as in the
original paper of Fleming and Viot \cite{Fleming:1979}, the genotype
of an individual (where genotype is a continuous instead of a discrete
variable). To our knowledge, the specific form of the limiting Fleming-Viot process above has not
previously been studied. In particular, although infinite dimensional
stochastic processes have been applied to multilevel population
dynamics of a single type \cite{Dawson:1991}, this appears to be the
first Fleming-Viot process for the evolution of two types under
opposing forces of selection at two population scales.

The dynamical properties of the deterministic partial differential equation \eqref{det_pde} are the focus of the next section. The proofs of weak convergence (Theorems \ref{det_thm} and \ref{fv_thm}) are deferred to section \ref{weak_convergence}.

\section{Properties of the deterministic limit}\label{det_steady_states}
We begin with a closed-form expression for solutions to the deterministic partial differential equation \eqref{det_pde}.
\begin{lemma}\label{method_char}
The solution to the deterministic partial differential equation \eqref{det_pde} with initial measure $\mu_0$ is given by
\begin{align}
\label{pde_soln} \mu_t(dx) & = (G_t\mu_0)(dx) = (\mu_0\phi^{-1}_t)(dx)\cdot w_t(x)
\end{align}
where 
\[\phi^{-1}_t(x) = \frac{x}{e^{-t}+x(1-e^{-t})}\]
\[w_t(x) = \left[(e^{-t}+x(1-e^{-t}))e^{t-\int_0^th(z)dz}\right]^{\lambda}\]
and $h(t)$ satisfies $h(t) = \inprod{x}{\mu_t}$
\end{lemma}
\textbf{Remark 1:} $(\mu_0\phi^{-1}_t)(dx):=\mu_0(\phi^{-1}_t(dx))$ captures the changes in the initial data that are solely due to the flux term. This expression is also known as the push-forward measure of $\mu_0$ under the dynamics of $\phi$. As we will see in the proof, $\phi_t(x)$ is precisely the characteristic curve for the spatial variable $x$ and includes a normalizing constant.
The multiplication by $w_t(x)$ captures the changes in the initial data that are due to the forcing term in \eqref{det_pde} and includes a normalizing factor. 

\textbf{Remark 2:} Density-valued solutions are given by
\begin{align}
\notag\eta_t(x) &= \eta_0(\phi_t^{-1}(x))\partial_x\phi_t^{-1}(x)\cdot w_t(x)\\
\label{density_soln}&=\eta_0\( \tfrac{x}{e^{-t}+x(1-e^{-t})}\)[e^{-t}+x(1-e^{-t})]^{\lambda-2}e^{(\lambda-1)t-\lambda\int_0^th(z)dz}
\end{align}
To see this, suppose $\mu_0(dx) = \eta_0(x)dx$. Then for any test function $f$,
\[\int_0^1 f(x) (\mu_0\phi^{-1}_t)(dx) = \int_0^1 (f\circ\phi_t)(x)\mu_0(dx) = \int_0^1 f(y)\eta_0(\phi_t^{-1}(y))\partial_y\phi_t^{-1}(y)dy\]
The first equality follows from the change-of-variable property of push-forward measures and the second from a standard change of variables. The limits of integration do not change because $0$ and $1$ are fixed points of both $\phi_t$ and $\phi_t^{-1}$.

\begin{proof}[Proof of Lemma~\ref{method_char}]
We apply the method of characteristics (see for example \cite{Pinchover:2005}) to obtain a formula for a density-valued solution. We then prove that the weak, measure-valued analog of this solution satisfies \eqref{det_pde}.
Consider the following modification of \eqref{nondim_strong_pde}:
\begin{align}
\label{h_det_pde}\fpart{}{t}\xi(t,x) = \fpart{}{x}[x(1-x)\xi(t,x)]+\lambda \xi(t,x)\left[x - h(t)\right]
\end{align}
where  $h(t)$ is a general function in time and $\xi_0\in C^1([0,1])$. Note that when $h(t) = \int_0^1 y\xi(t,y)dy$, this differential equation is equivalent to \eqref{nondim_strong_pde}. To be clear about which equation we are solving, we use $\xi(t,x)$ to denote solutions when $h(t)$ is unspecified.

Rewriting \eqref{h_det_pde}:\begin{align*}
\(\fpart{}{t}\xi, \fpart{}{x}\xi,-1\) \cdot \(1, -x(1-x),  \left[(1-2x) + \lambda(x-h(t))\right]\xi\)& =0
\end{align*}
The second vector is therefore tangent to the solution surface and gives the rates of change for the $t$, $x$, and $\xi$ coordinates. Let the initial condition be parameterized as $(0,x,\xi_0(x)) = (0,p,\xi_0(p))$. The $t$, $x$, and $\xi$ coordinates change according to the characteristic equations
\begin{align*}
\tfrac{dt}{dq} &= 1&t(0,p) = 0\\
\tfrac{dx}{dq}&=-x(1-x)&x(0,p) = p\\
\tfrac{d\xi}{dq}& = \left[(1-2x(q,p)) + \lambda(x(q,p)-h(t(q,p)))\right]\xi&\xi(0,p) = \xi_0(p)
\end{align*}
where $q$ is the parameter as we move through the solutions in time.
The first two ordinary differential equations have solutions
\begin{align}
\notag t(q,p) &= q\\
\label{phi_def}x(q,p) & = \frac{p}{p-(p-1)e^{q}}=:\phi_q(p)
\end{align}
From this, the third differential equation can be solved exactly:
\begin{align*}
\frac{d\xi}{dq}& = \left[1+\tfrac{p}{p-(p-1)e^{q}}(\lambda-2)- \lambda h(q)\right]\xi\\
\xi(q,p) & = \xi_0(p)\exp\left\{q -\lambda \int_0^q h(z)dz +(\lambda-2)\int_0^q \frac{p}{p-(p-1)e^{z}}~dz\right\}\\
 & = \xi_0(p)e^{q -\lambda \int_0^q h(z)dz} \left[p(e^{-q}-1)+1\right]^{(-\lambda+2)}
\end{align*}
Next, make the substitutions $q=t$ and $p = \phi_t^{-1}(x)$ from \eqref{phi_def} to obtain $\xi$ in terms of $t$ and $x$:
\begin{align}
\notag\xi(t,x) 
& = (\xi_0\circ\phi_t^{-1})(x) \left[e^{-t}+x(1-e^{-t})\right]^{(\lambda-2)}e^{(\lambda-1)t -\lambda \int_0^t h(z)dz}\\
\label{density_soln_xi}& =  (\xi_0\circ\phi_t^{-1})(x) \partial_x\phi_t^{-1}(x)\cdot w_t(x)
\end{align}

If $h(t)$ satisfies $h(t) = \int_0^1y\xi(t,y)dy$, then by definition, $\xi(t,x)$ solves the partial differential equation \eqref{nondim_strong_pde}. Conversely, if $\xi(t,x)$ solves the partial differential equation \eqref{nondim_strong_pde}, it also solves the differential equation \eqref{h_det_pde} with $h(t) = \int_0^1y\xi(t,y)dy$. Therefore this above expression, along with the condition $h(t)=\int_0^1 y\xi(t,y)dy$, are equivalent to solutions of \eqref{nondim_strong_pde}.

To extend this result to measures, suppose we have a strong solution $\eta_t(x)$ with initial condition $\eta_0$:
\[\eta_t(x) = (\eta_0\circ\phi_t^{-1})(x) \partial_x\phi_t^{-1}(x)\cdot w_t(x)\]
Using a similar calculation as that in Remark~2, the measure $\mu_t(dx)$ corresponding to $\eta_t(x)$ is given by 
\[\mu_t(dx) = (\mu_0\phi_t^{-1})(dx)\cdot w_t(x)\]
It remains to check that this satisfies the weak deterministic partial differential equation \eqref{det_pde}, with $h(t)= \inprod{x}{\mu_t}$. The left hand side of the equation is
\begin{align*}
\tfrac{d}{dt}\inprod{f}{\mu_t}& = \tfrac{d}{dt}\int_0^1 f(x)w_t(x)(\mu_0\phi_t^{-1})(dx) = \tfrac{d}{dt}\int_0^1f(\phi_t(x))w_t(\phi_t(x))\mu_0(dx)
\end{align*}
Differentiating under the integral sign, expanding out the expressions for $\partial_t\phi_t$ and $\partial_t(w_t(\phi_t(x))$, and applying change of variables for push-forward measures again, we obtain
\begin{align*}
\tfrac{d}{dt}\inprod{f}{\mu_t}& = {- \int_0^1 x(1-x) f'(x)w_t(x)(\mu_0\phi_t^{-1})(dx)} + \lambda\int_0^1\left[x-h(t)\right]f(x)w_t(x)(\mu_0\phi_t^{-1})(dx)
\end{align*}
This matches right hand side of the weak deterministic partial differential equation \eqref{det_pde}.
\end{proof}

In practice, the condition $h(t) =\inprod{x}{\mu_t}$ is difficult to use. The following provides an equivalent and simpler condition.

\begin{lemma}[Conservation of measure condition]
\label{cons_cond}
Suppose $\xi$ is a weak measure-valued solution to the deterministic partial differential equation \eqref{h_det_pde} with initial condition $\int_0^1\xi_0( dx) = 1$. Then
\[h(t) = \int_0^1 y\xi(t,dy)\qquad \text{ if and only if }\qquad \int_0^1\xi(t,dy)=1~\forall~t>0\]
\end{lemma}
\begin{proof}
($\Rightarrow$ direction) Suppose $h(t) = \int_0^1y\xi(t,dy)$. Then $\xi$ is a weak measure-valued solution to \eqref{det_pde}. Taking the test function $f\equiv1$, we obtain
\begin{align*}
\tfrac{d}{dt}\inprod{1}{\xi} & = 0+\lambda[\inprod{x}{\xi}-\inprod{1}{\xi}\inprod{x}{\xi}]=0
\end{align*}
Thus, if the initial data has total measure $1$, $\inprod{1}{\xi}$ remains constant at $1$ for all $t\geq0$.

($\Leftarrow$ direction) Suppose $\int_0^1\xi(t,dx)=1$ for all $t>0$. Again take the test function $f\equiv 1$ but this time with unspecified $h(t)$:
\begin{align*}
0=\tfrac{d}{dt}\inprod{1}{\xi} &= 0+\lambda[\inprod{x}{\xi} -\inprod{1}{\xi}h(t)]=\lambda[\inprod{x}{\xi}-h(t)].
\end{align*}
For this to hold, we must have $h(t) = \int_0^1 x\xi(t,dx)$. 
\end{proof}

The above lemmas imply that solutions $\mu_t(dx)$ to \eqref{det_pde}
can be obtained by using formula \eqref{pde_soln} from
Lemma~\ref{method_char} and imposing the conservation of measure
condition $\inprod{1}{\mu_t}\equiv1$ from Lemma~\ref{cons_cond}. We
illustrate this with some examples of exactly
solvable solutions for special choices of initial data. We will see
that the  long time behavior of the  examples is consistent with results stated in Theorem~\ref{thm:steadyState}.

\begin{example}\label{ex1}
Initial measure concentrated at $x_0\in[0,1]$, i.e.~$ \mu_0=\delta_{x_0}$\\
Using formula \eqref{pde_soln},
\begin{align*}
\int f(x)\mu_t(dx) &= \int f(x) w_t(x)\delta_{x_0}(\phi_t^{-1}(dx))= f(\phi_t(x_0))w_t(\phi_t(x_0))\\
&=\int f(x)w_t(x)\delta_{\phi_t(x_0)}(dx)
\end{align*}
Thus $\mu_t(dx) = w_t(x)\delta_{\phi_t(x_0)}(dx)$. Imposing the conservation of measure condition gives $\mu_t(dx) = \delta_{\phi_t(x_0)}(dx)$. In other words, an initial delta measure at $x_0$ moves as a delta measure along the $x$ axis with position given by $\phi_t(x_0)$, the solution to the negative logistic equation with initial position $x_0$.
\end{example}

\begin{example}
Initial uniform density: $\eta_0(x) = 1$, i.e~$\mu_0(dx) = dx$\\
Using formula \eqref{density_soln},
\begin{align*}
\eta_t(x) & = e^{(\lambda-1)t -\lambda \int_0^t h(z)dz} \left[e^{-t}+x(1-e^{-t})\right]^{(\lambda-2)}
\end{align*}
Imposing conservation of measure:
\begin{align*}
e^{(\lambda-1)t -\lambda \int_0^t h(z)dz} &= \left[\int_0^1 \left[e^{-t}+x(1-e^{-t})\right]^{(\lambda-2)} dx\right]^{-1}\\
& =
\begin{cases}\displaystyle
\frac{(\lambda-1)(1-e^{-t})}{1-e^{-(\lambda-1)t}} & \text{if }
\lambda\neq1\\
\\
\displaystyle\frac{1-e^{-t}}{t} & \text{if } \lambda=1   
\end{cases}
\end{align*}
Thus,
\begin{align*}
\eta_t(x) & =
\begin{cases}\displaystyle
\frac{(\lambda-1)(1-e^{-t})}{1-e^{-(\lambda-1)t}} \left[e^{-t}+x(1-e^{-t})\right]^{(\lambda-2)}
 & \text{if } \lambda\neq1\\
\\
\displaystyle \frac{1-e^{-t}}{t}\left[e^{-t}+x(1-e^{-t})\right]^{(\lambda-2)} & \text{if } \lambda=1
\end{cases}
\end{align*}

Note that $\eta_0\equiv 1$ corresponds to an initial condition
satisfying the hypothesis of Theorem~\ref{thm:steadyState} with
$\alpha=1$. As predicted when $\lambda>1$, we obtain
$\eta(t,x)\to(\lambda-1)x^{\lambda-2} = \text{Beta}(\lambda-1,
1)$ as  $t\to\infty$.

\label{ex2}
\end{example}
 The following is an example with 
$\alpha>1$. 
\begin{example}\label{ex3}
If  $\eta_0(x) = 2(1-x)$, i.e.~$\mu_0([1-x,1])=x^2$, then the corresponding $\alpha$ from
Theorem~\ref{thm:steadyState}  is  $\alpha=2$.\\

Using formula~\eqref{density_soln}
\begin{align*}
\eta_t(x) & = 2e^{(\lambda-2)t -\lambda \int_0^t h(z)dz} (1-x)\left[e^{-t}+x(1-e^{-t})\right]^{(\lambda-3)}
\end{align*}
Imposing the condition in Lemma~\ref{cons_cond} to solve for the $h(z)$ term
\begin{align*}
  e^{(\lambda-2)t -\lambda \int_0^t h(z)dz} = &\left[2\int_0^1 (1-x)\left[e^{-t}+x(1-e^{-t})\right]^{(\lambda-3)} dx\right]^{-1}\\
  = &
 \begin{cases}
   \frac{(\lambda-2)(1-e^{-t})}{2}\left[\frac{1}{(\lambda-1)(1-e^{-t})} - \frac{e^{-(\lambda-1)t}}{(\lambda-1)(1-e^{-t})}-e^{-(\lambda-2)t}\right]^{-1} & \text{if } \lambda\neq2\\
   \frac{(1-e^{-t})^2}{2te^{-t}} & \text{if } \lambda=2
 \end{cases}
\end{align*}

As predicted by Theorem~\ref{thm:steadyState} for $\lambda>2=\alpha$,

\begin{align*}
  \eta_t(x) \to \tfrac{1}{2}(\lambda-2)(\lambda-1)(1-x)x^{\lambda-3} =
  \text{Beta}(\lambda-2, 2)
\end{align*}
as $t\to\infty$. 
\end{example}

\begin{example}\label{ex4}
$\eta_0 (x)=\frac{1}{c}\cdot 1_{[0,c]}(x)$ with $c<1$. 

Using formula~\eqref{density_soln}
\[\eta_t(x)=\tfrac{1}{c}\mathbf{1}\{x\leq \phi_t(c)\}w_t(x)\partial_x\phi_t^{-1}(x)
\]
Since $\phi_t(c) = \frac{ce^{-t}}{1-c+ce^{-t}}\to 0$ as $t\to\infty$,
$\eta_t(x)\to 0$ for any $x>0$. Since $\eta$ must have total mass 1, it
follows that regardless of the value of $\lambda$,  $\eta_t(x)dx\to
\delta_0(dx) \text{ for any } c<1$. This can also be seen by applying Theorem~\ref{thm:steadyState}
and noting that  $\mu_0([1-c,1])=\int_{1-c}^1\eta_0(x)dx = 0$.

\end{example}

We end these examples with solutions for $\mu_0$ that are mixtures of
delta measures and densities. First, note that it is straightforward
to extend Example~\ref{ex1} to the case where $\mu_0(dx) = \sum a_i
\delta_{x_i}(dx)$ is a linear combination of delta measures, $a_i>0$
for all $i$. Applying \eqref{pde_soln}, we obtain
\[\mu(t,dx) = \sum_i a_iw_t(x)\delta_{\phi_t(x_i)}(dx) = \sum_i a_i(t)\delta_{x_i(t)}(dx)
\]
where $x_i(t) = \phi_t(x_i)$ and $a_i(t) = a_i w_t(x)|_{x=x_i(t)}$. Our earlier system of equations \eqref{eq:particles} is obtained from this and the definitions of $\phi_t(x)$ and $w_t(x)$.

Second, we consider a combination of a delta measure and a density
\[\mu_0(dx) = a\delta_{x_0}(dx) + (1-a) v_0(x)dx\]
Notice that the formula for the solution \eqref{pde_soln} at first seems linear in the initial condition:
\begin{align*}
\int f(x) \mu_t(dx) &= \int f(x) (G_t\mu_0)(dx)\\
& = \int f(x) w_t(x) [a\delta_{\phi_t(x_0)}(dx) + (1-a) v_0(\phi^{-1}_t(x))\partial_x\phi^{-1}_t(x)dx]\\
& = \int f(x) [a(G_t \delta_{x_0})(dx) + (1-a)(G_t v_0)(dx)]
\end{align*}
This gives $(G_t \mu_0)(dx) = a(G_t\delta_{x_0})(dx) + (1-a)(G_tv_0)(dx) $. However, this notation is misleading because implicit in the $G_t$ operator is the function $h(t)$, the mean of the overall process over time. Here, $h(t)$ involves both  the delta measure and the density. The solution operator $G_t$ is therefore not linear for this reason.

Nevertheless, we can still use this formula to obtain expressions for solutions. We illustrate this with a concrete example.

\begin{example}\label{ex5}
  Take $x_0=0$ and $v_0(x)$ the density function for  $\text{Beta}(\lambda-\alpha, \alpha)$ with $ \alpha \in (0,\lambda)$. Using the solution formula and direct calculation, we obtain
\begin{align*}
\mu_t(dx) &= aw_t(0)\delta_0(dx)  + (1-a)w_t(x)(v_0\circ\phi_t^{-1})(x)\partial_x\phi_t^{-1}(x)dx\\
& = e^{-\lambda\int_0^t h(z)dz)}\left\{a \delta_0(dx) + (1-a)e^{(\lambda-\alpha)t}v_0(x)dx\right\}
\end{align*}
Note in particular that $\mu_t$ remains a linear combination of $\delta_0$ and the Beta distribution. The Beta distribution ultimately dominates because $\lambda>\alpha$.
\end{example}

We now use Lemma~\ref{method_char} to show that Beta distributions,
$\delta_0$, and $\delta_1$ are fixed points for the deterministic
partial differential equation and thus provide a proof of
Lemma~\ref{lem:fixedPoints} announced earlier in this note.
\begin{proof}[Proof of Lemma~\ref{lem:fixedPoints}]  Note that we could prove this lemma by substituting $\delta_0$, $\delta_1$, and the Beta distribution into the deterministic partial differential equation \eqref{det_pde} and showing the right-hand side equals zero. Instead, we will show that these distribution are fixed points of the solution operator. Let $v$ be the density of the Beta distribution, \[v(x) = \frac{1}{B(\lambda-\alpha, \alpha)}x^{\lambda-\alpha-1}(1-x)^{\alpha-1}.\] The mean of $v$ is $\frac{\lambda-\alpha}{\lambda}$. Using \eqref{density_soln}

\begin{align*}
(G_tv)(x) &= v\( \tfrac{x}{e^{-t}+x(1-e^{-t})}\)[e^{-t}+x(1-e^{-t})]^{\lambda-2}e^{(\lambda-1)t-(\lambda-\alpha)t}= v(x)
\end{align*}
$v$ is therefore a fixed point of the solution operator and hence is a fixed point of the deterministic partial differential equation.

For $\delta_0$ and $\delta_1$, we use Example~\ref{ex1} above to obtain $(G_t\delta_{x_0})(dx) = \delta_{\phi_t(x_0)}(dx)$. Since $x_0=0$ and $x_0=1$ are fixed points of $\phi_t$, it follows that $\delta_0$ and $\delta_1$ are fixed points of $G_t$. 
\end{proof}

We now prove when
the fixed points are stable.
We begin with a lemma which gives more general conditions than those given in Theorem~\ref{thm:steadyState} for the delta
measure at zero to attract a given initial condition.
\begin{lemma}\label{lem:delta0Attracts}If for some $\alpha\geq
  \lambda>0$,
  \begin{align*}
  \lim_{x \rightarrow 0}  x^{-\alpha}  \mu_0([1-x,1]) < \infty
  \end{align*}
then $\mu_t \rightarrow \delta_0$ as $t \rightarrow \infty$. In
particular, this condition holds if $\mu_0([1-\epsilon,1])=0$ for some $\epsilon>0$.
\end{lemma}

To prove this and subsequent results, we will need the following
technical lemma.
\begin{lemma}\label{lem:hTech}
  Setting $h(t)=\inprod{x}{\mu_t}$, the following two implications
  hold:
  \begin{align*}
    \int_0^\infty h(t)\,dt < \infty \qquad&\Longrightarrow \qquad h(t) \rightarrow0 
    \quad\text{as}\quad t\rightarrow \infty\,.\\
    \int_0^\infty \big[1-h(t)\big]\,dt < \infty \qquad&\Longrightarrow \qquad h(t) \rightarrow1 
    \quad\text{as}\quad t\rightarrow \infty\,.
  \end{align*}
\end{lemma}
\begin{proof}[Proof of Lemma~\ref{lem:hTech}]
  Since $h(t) \geq 0$ and $1-h(t) \geq 0$, the only obstruction to the
  implication is that $h(t)$ (or $1-h(t)$) could have ever shorter and
  shorter intervals were they return to an order one value before
  returning to a value close to zero. This would require $h(t)$ to
  have unbounded derivatives. However this is not possible since
  \begin{align*}
    \frac{d h}{dt}(t) =  -\big(h - \inprod{x^2}{\mu_t}\big) + \lambda \big(
    \inprod{x^2}{\mu_t} - h^2 \big)
  \end{align*}
from which one easily see that $ -1 \leq \frac{d h}{dt}(t)  \leq \lambda$
since $0 \leq h - \inprod{x^2}{\mu_t}\leq 1$ and $0\leq \inprod{x^2}{\mu_t}- h^2\leq 1$.
\end{proof}

\begin{proof}[Proof of Lemma~\ref{lem:delta0Attracts}]
As usual let $h(t)=\inprod{x}{\mu_t}$. We begin by observing that if
\begin{align*}
  \int_0^\infty h(t) dt < \infty 
\end{align*}
then $h(t) \rightarrow 0$ as $t \rightarrow \infty$ by Lemma~\ref{lem:hTech} and $\mu_t
\rightarrow \delta_0$ as we wish to prove. Thus, we henceforth assume
that  $\int_0^\infty h(t) dt =\infty$. 
Under this assumption, we will show that for any continuous function $f$
\begin{align*}
  \int_0^1 f(x) \mu_t(dx)  \rightarrow f(0)\qquad\text{as}\qquad t
  \rightarrow \infty\,.
\end{align*}
 Since $f$ is continuous, given any $\epsilon>0$, there exists a
 $\delta>0$ so that $|f(x)-f(0)|<\epsilon$ whenever $x\leq
 \delta$. Hence
 \begin{equation}\label{epsBound}
  \Big| \int_0^1 f(x) \mu_t(dx) -f(0)\Big| \leq     \int_0^1 |f(x)
  -f(0)| \mu_t(dx) \leq \epsilon + \int_\delta^1 |f(x)
  -f(0)| \mu_t(dx)
 \end{equation}
Now setting 
\begin{align*}
  \int_\delta^1 |f(x)
  -f(0)| \mu_t(dx) &=   \int_{\phi^{-1}_t(\delta)}^1  |(f\circ \phi_t)(x)
  -f(0)|\,(w_t\circ \phi_t)(x)\, \mu_0(dx)\\
  &\leq 2 \|f\|_\infty \int_{\phi_t^{-1}(\delta)}^1  ( w_t\circ \phi_t)(y)\,
  \mu_0(dy).\\
\end{align*}
Since for all $y \in [\phi_t^{-1}(\delta),1]$ and $t > 0$, we have
\begin{align*}
   ( w_t\circ \phi_t)(y) \leq
   e^{\lambda t-\lambda \int_0^t h(s) ds}
\end{align*}
we see that 
\begin{align*}
   \int_\delta^1 |f(x)
  -f(0)| \mu_t(dx) & \leq  2 \|f\|_\infty  
  e^{\lambda  t -\lambda \int_0^t h(s) ds} \mu_0([\phi_t^{-1}(\delta),1])\,.
\end{align*}
Now using the assumptions on $\mu_0$ and that $\phi_t^{-1}(\delta)
\geq 1-D e^{-t}$ for some $D>0$ and all $t>0$, one has that
\begin{align*}
  e^{\lambda  t -\lambda \int_0^t h(s) ds} \mu_0([\phi_t^{-1}(\delta),1]) \leq \hat D
  e^{-(\alpha-\lambda) t -\lambda \int_0^t h(s) ds}
\end{align*}
for some constant $\hat D$ and all $t>0$. 
Since $\alpha\geq \lambda$ and $\int_0^\infty h(s) ds = \infty$, this
bound converges to zero as $t \rightarrow \infty$ and the proof is
complete as the $\epsilon$ in \eqref{epsBound} was arbitrary.
\end{proof}

\begin{proof}[Proof of Theorem~\ref{thm:steadyState}] We start with
  the setting when $\mu_0(\{1\})>0$ and begin by writing $\mu_t(dx) =
  a_t \delta_1(dx) + (1-a_t) \nu_t(dx)$ for some time dependent process
  $a_t \in [0,1]$ with $a_0 >0$ and some probability measure valued
  process $\nu_t(dx)$. As usual we define $h(t) = \inprod{x}{\mu_t}$ and using the representation given in \eqref{pde_soln},
  one sees that $a_t$ solves
  \begin{equation*}
    \frac{d a_t}{dt} = \lambda a_t\big( 1 - h(t)\big) 
    \quad\Longrightarrow\quad a_t = a_0 \exp\Big( \lambda\int_0^t [1-h(s)]ds \Big). 
  \end{equation*}
Since $1 - h(t)\geq 0$, we know that $\int_0^t [1-h(s)]ds$ converges 
as $t \rightarrow \infty$. If it converges to $\infty$ then $a_t$ also 
converges to $\infty$ since $a_0 >0$. However this is impossible since $a_t \in 
[0,1]$ for all $t\geq 0$. Thus, we conclude that  $\int_0^t 
[1-h(s)]ds<\infty$. Then Lemma~\ref{lem:hTech} implies that 
$h(t)\rightarrow 1$ which in turn implies that $\mu_t \rightarrow \delta_1$ as 
$t \rightarrow \infty$. 

We know turn to the setting when $x^{-\alpha} \mu_0([1-x,1])
\rightarrow C>0$ as $x \rightarrow 0$. The case when $\lambda \leq\alpha$ is already handled by
Lemma~\ref{lem:delta0Attracts} leaving only  the case when
$\lambda>\alpha>0$ to be proven.
For $x \in [0,1]$, define $U(x)=\mu_0([0,x])$. Since $\mu_0$ is a
probability measure we know that $U$ has finite variation and is
regular in the sense that both the right limit $U(x^+)$ and the left
limit $U(x^-)$ exist, where $U(x^\pm) = \lim U(y)$ as $y
\rightarrow^\pm x$. At the extreme points, only the limit obtained by staying in $[0,1]$ 
is defined.

Now for any smooth function $f$ of $[0,1]$, we have from \eqref{pde_soln} that 
\begin{align*}
  \int_0^1 f(x) \mu_t(dx) = Z_t \int_0^1 f(x)g_t(x) 
  (\mu_0\phi^{-1}_t)(dx)= Z_t \int_0^1 [(fg_t)\circ\phi_t](x)\,
  \mu_0(dx) 
\end{align*}
where $w_t(x)$ has been written as the product of $g_t(x)=(e^{-t}+x(1-e^{-t}))^\lambda$ and $Z_t$ some positive,
time dependent normalizing constant. It is enough to show that for 
some time positive, dependent constant $K_t$,
\begin{equation}\label{eq:enoughToShow}
  K_t  \int_0^1 [(fg_t)\circ\phi_t](x) 
  \mu_0(dx) \rightarrow \int_0^1 
  f(x) x^{\lambda-\alpha-1}(1-x)^{\alpha-1}dx \quad\text{as}\quad t 
  \rightarrow \infty\,. 
\end{equation}

Since $x \mapsto f(x) g_t(x)$ is continuous on $[0,1]$, even if $U(x)$
has discontinuities the integration by parts formula for 
Lebesgue-Stieltjes integrals produces 
\begin{align*}
  \int_0^1 [(fg_t)\circ\phi_t](x)\,  \mu_0(dx)&= (fg_tU)(1^-) - (fg_tU)(0^+) -
  \int_0^1 \partial_x[(fg_t)\circ\phi_t](x)\, U(x)\, dx\\
&= [fg_t(U-1)](1^-) + [fg_t(1-U)](0^+) +
  \int_0^1 \partial_x[(fg_t)\circ\phi_t](x)\, [1-U](x)\, dx\,. 
\end{align*}
Here we have used that  $\phi_t$ is continuous with $\phi_t(1)=1$ and  $\phi_t(0)=0$. 

 First observe that $1-U(1^-)=0$ since $\mu_0([1-x,1])\rightarrow 0$ as $x\rightarrow 0$ by assumption and that $g_t(0^+)= e^{-\lambda t}$. Hence 
 \begin{equation}
   \label{eq:boundaryTerms}
  [fg_t(1-U)](1^-) + [fg_t(U-1)](0^+)= [U(0^+)-1] f(0) e^{-\lambda t}.
 \end{equation}
Now turning to the integral term, applying the chain rule and 
changing variables to $y=\phi_t(x)$ produces 
\begin{align*}
    \int_0^1 \partial_x[(fg_t)\circ\phi_t](x)\, [1-U](x)\, dx &=   \int_0^1 [\partial_x(fg_t)\circ\phi_t](x)\, [1-U](x)\,(\partial_x \phi_t)(x)\,dx\\&=
    \int_0^1[\partial_x(fg_t)](y) 
    \,[(1-U)\circ  \phi^{-1}_t](y) dy 
\end{align*}
For any fixed $x \in (0,1)$ by 
direct calculation and use of the assumption on $\mu_0$, one sees that 
\begin{equation*}
  \left.\begin{aligned}
\partial_x(fg_t)(x) &\rightarrow \partial_x(x^\lambda f)(x)\\
 e^{\alpha t}( 1- U(\phi_t^{-1}(x)) )= e^{\alpha t}\mu_0([\phi_t^{-1}(x),1])&\rightarrow  C\big(\tfrac{1-x}{x}\big)^{\alpha}
  \end{aligned}\right\}
 \quad\text{as}\quad t \rightarrow \infty\,. 
\end{equation*}
Combining these facts with \eqref{eq:boundaryTerms} and the fact that $e^{-(\lambda-\alpha)t} \rightarrow 0$ as $t \rightarrow \infty$ since $\lambda >\alpha$ produces
\begin{align*}
   e^{\alpha t} \int_0^1 [(fg_t)\circ\phi_t](x)\,  \mu_0(dx) \rightarrow 
   C\int_0^1 \partial_x(x^{\lambda} f)(x) 
   \Big(\frac{1-x}{x}\Big)^{\alpha} dx \quad\text{as}\quad t\rightarrow \infty\,. 
\end{align*}
for some new positive constant $C$. 
Now since integration by parts implies that 
\begin{align*}
 \frac1\alpha \int_0^1 \partial_x(x^{\lambda} f)(x) 
   \Big(\frac{1-x}{x}\Big)^{\alpha} dx= \int_0^1f(x)  x^{\lambda-\alpha-1}(1-x)^{\alpha-1} dx 
\end{align*}
the last part of the proof is complete.
\end{proof}

\section{Proofs of weak convergence}\label{weak_convergence}
The proofs of Theorems \ref{det_thm} and \ref{fv_thm} follow a standard procedure \cite{Joffe:1986,Fournier:2004, Champagnat:2006}. Both proofs require: (i) tightness of the sequence of stochastic processes -- which implies a subsequential limit, and (ii) uniqueness of this limit. For the tightness of $\{\mu^{m,n}_t\}_{m,n}$ on $D([0,T], \mathcal{P}([0,1]))$, it is sufficient, by Theorem~14.26 in Kallenberg \cite{Kallenberg:1997} to show that $\{\inprod{f}{\mu^{m,n}_t}\}$ is tight on $D([0,T], \mathbb{R})$ for any test function $f$ from a countably dense subset of continuous, positive functions on $[0,1]$. For the uniqueness of solutions to the partial differential equation in Theorem~\ref{det_thm}, we apply Gronwall's inequality. For uniqueness of solutions to the martingale problem in Theorem~\ref{fv_thm}, we apply a Girsanov theorem by Dawson \cite{Dawson:2010}.

\subsection{Semimartingale property of multilevel selection process}
It will be useful for what follows to treat $\inprod{f}{\mu^{m,n}_t}$ as a semimartingale. Below, $D^+_xf$ is the first order difference quotient of $f$ taken from the right,  $D^-_xf$ is the first order difference quotient of $f$ taken from the left, and $D_{xx}f$ is the second order difference quotient.
\begin{lemma}
\label{semi_mart_lem}
For $f\in C^2([0,1])$ and $\mu^{m,n}_t$ with generator $L^{m,n}$ defined in \eqref{gen_eqn},
\begin{align}
\label{semimart} \inprod{f}{\mu^{m,n}_t}-\inprod{f}{\mu^{m,n}_0}& =  A^{m,n}_t(f) + M^{m,n}_t(f)
\end{align}
where $A^{m,n}_t(f)$ is a process of finite variation, $A^{m,n}_t(f): = \int_0^t a^{m,n}_z(f)dz$, with
\begin{align}
 \label{finite_var}a^{m,n}_t(f) & =  \sum_i \mu^{m,n}_t(\tfrac{i}{n})\tfrac{i}{n}(1-\tfrac{i}{n})\left[\tfrac{1}{n}D_{xx}f(\tfrac{i}{n}) -sD^-_{x}f(\tfrac{i}{n})\right]\\
\notag&\qquad+ wr\left\{\sum_{j}\mu^{m,n}_t(\tfrac{j}{n})\tfrac{j}{n}f(\tfrac{j}{n})-\sum_{i}\mu^{m,n}_t(\tfrac{i}{n})f(\tfrac{i}{n})\sum_j\mu^{m,n}_t(\tfrac{j}{n})\tfrac{j}{n}\right\}
\end{align}
and $M^{m,n}_t(f)$ is a c\`adl\`ag  martingale with (conditional) quadratic variation
\begin{align}
\notag\langle M^{m,n}(f)\rangle_t =& \tfrac{1}{m}\int_0^t \left\{\tfrac{1}{n}\sum_i \mu^{m,n}_z(\tfrac{i}{n})\tfrac{i}{n}(1-\tfrac{i}{n})\left[\(D^+_xf(\tfrac{i}{n})\)^2+(1+s)\(D^-_xf(\tfrac{i}{n})\)^2\right]\right.\\
\label{quad_var}&\qquad \left.+w\sum_{i,j}\mu^{m,n}_z(\tfrac{i}{n})\mu^{m,n}_z(\tfrac{j}{n})(1+r\tfrac{j}{n})(f(\tfrac{i}{n})-f(\tfrac{j}{n}))^2\right\} dz
\end{align}
\end{lemma}

\begin{proof}
By Dynkin's formula (see, for example, Lemma~17.21 in \cite{Kallenberg:1997}), 
\[\psi(\mu^{m,n}_t)-\psi(\mu^{m,n}_0) - \int_0^t (L^{m,n}\psi)(\mu^{m,n}_s)ds\]
where $\psi\in dom(L^{m,n})$, is a c\`adl\`ag martingale. In particular, this is true for \[\psi(\mu^{m,n}_t)=F(\inprod{f}{\mu^{m,n}_t})\]
where $f\in C^2([0,1])$ and $F:\mathbb{R}\to\mathbb{R}$.
Setting $F(x) = x$ and plugging this $f$ into \eqref{gen_eqn}:
\begin{align*}
(L^{m,n}\inprod{f}{\cdot})(v)
= &\sum_i v(\tfrac{i}{n})\tfrac{i}{n}(1-\tfrac{i}{n})\left[\tfrac{1}{n}D_{xx}f(\tfrac{i}{n}) -sD^-_{x}f(\tfrac{i}{n})\right]\\
&\qquad+ wr\left\{\sum_{j}v(\tfrac{j}{n})\tfrac{j}{n}f(\tfrac{j}{n})-\sum_{i}v(\tfrac{i}{n})f(\tfrac{i}{n})\sum_jv(\tfrac{j}{n})\tfrac{j}{n}\right\}
\end{align*}
Thus,
\begin{align}
\label{mart}\inprod{f}{\mu^{m,n}_t}-\inprod{f}{\mu^{m,n}_0} - \int_0^t a^{m,n}_z(f) dz = M^{m,n}_t(f)
\end{align}
 where $M^{m,n}_t(f)$ is some martingale and $a^{m,n}_t (f) = (L^{m,n}\inprod{f}{\cdot})(\mu^{m,n}_t)$. $A_t(f)$ is a process of finite variation because for a given $f$, $a_t^{m,n}(f)$ is uniformly bounded in $t$.

Next, setting $F(x) = x^2$ and plugging this $\psi$ into \eqref{gen_eqn}:  
\begin{align*}
(L^{m,n}\inprod{f}{\cdot}^2)(v)
 = &  2\inprod{f}{v}a^{m,n}_t(f)+\tfrac{1}{mn}\sum_i v(\tfrac{i}{n})\tfrac{i}{n}(1-\tfrac{i}{n})\left[\(D^+_xf(\tfrac{i}{n})\)^2+(1+s)\(D^-_xf(\tfrac{i}{n})\)^2\right]\\
&\qquad+\tfrac{w}{m}\sum_{i,j}v(\tfrac{i}{n})v(\tfrac{j}{n})(1+r\tfrac{j}{n})(f(\tfrac{i}{n})-f(\tfrac{j}{n}))^2
\end{align*}

Thus,
\begin{align}
\label{martqv_1}\inprod{f}{\mu^{m,n}_t}^2-\inprod{f}{\mu^{m,n}_0}^2 - \int_0^t c^{m,n}_z(f) dz=\text{martingale}
\end{align} 
where $c^{m,n}_t(f) = (L^{m,n}\inprod{f}{\cdot}^2)(\mu^{m,n}_t)$.

Alternatively, take $Y_t = \inprod{f}{\mu^{m,n}_t}$ and apply Ito's formula (for example, p78 in \cite{Protter:2004}) to $Y_t^2$ to obtain
\begin{align}
\label{martqv_2} \inprod{f}{\mu^{m,n}_t}^2-\inprod{f}{\mu^{m,n}_0}^2 =& 2\int_0^t\inprod{f}{\mu_z}a^{m,n}_z (f)dz + [M^{m,n}(f)]_t + \text{martingale}
\end{align}
where $[M^{m,n}(f)]_t$ is the quadratic variation process of $M^{m,n}_t$. Since  $\langle M^{m,n}(f)\rangle_t$ is the compensator of $[M^{m,n}(f)]_t$,  
\[[M^{m,n}(f)]_t-\langle M^{m,n}(f)\rangle_t\]
is a martingale. Thus,
\begin{align}
\label{martqv_3}\inprod{f}{\mu^{m,n}_t}^2-\inprod{f}{\mu^{m,n}_0}^2 - 2\int_0^t\inprod{f}{\mu^{m,n}_z}a^{m,n}_z (f)dz - \langle M^{m,n}(f)\rangle_t =& \text{ martingale}
\end{align}
The compensator $\langle M^{m,n}(f) \rangle_t$ is a predictable process of finite variation (see p118 in\cite{Protter:2004}). By the Doob-Meyer inequality (p103 in \cite{Protter:2004}), the martingale in \eqref{martqv_3}  is the same as the martingale in \eqref{martqv_1} .  Equating these martingale parts we obtain
\begin{align}
\label{equate_mart}  2\int_0^t\inprod{f}{\mu^{m,n}_z}a^{m,n}_z (f)dz + \langle M^{m,n}(f)\rangle_t&=\int_0^t c^{m,n}_z(f)dz\,.
\end{align}
Substituting in the expressions for $a^{m,n}_z$ and $c^{m,n}_z$ then gives the explicit expression for the conditional quadratic variation \eqref{quad_var} in the statement of the lemma.
\end{proof}

\subsection{Proof of deterministic limit}\label{det_limit}

To prove Theorem \ref{det_thm}, we need the two following lemmas. The first uses criteria in Billingsley \cite{Billingsley:1999} to show tightness of the sequence of processes $\inprod{f}{\mu^{m,n}_t}$. The second uses Gronwall's inequality to show uniqueness of solutions to the limiting system.

\begin{lemma}\label{det_tightness}
The processes $\inprod{f}{\mu^{m,n}_t}$, as a sequence in $\{(m,n)\}$, is tight for all positive-valued test functions $f\in C^1([0,1])$.
\end{lemma}
\begin{proof}
By Theorem~13.2 in \cite{Billingsley:1999}, a sequence of probability measures $\{P_n\}$ on $D([0,T], \mathbb{R}^+)$ is tight if and only if
 (i) for all $\eta>0$, there exists $a$ such that
\[P_{n}\(x:\sup_{t\in[0,T]}|x(t)|\geq a\) \leq \eta\text{  for }n\geq 1\] and (ii) for all $\epsilon>0$ and $\eta>0$, there exists $\delta\in(0,1)$ and $n_0$ such that
\[P_{n}(x: w_x'(\delta)\geq \epsilon)\leq\eta\text{  for all  }n>n_0\]
where $w'$ is the modulus of continuity for c\`adl\`ag processes and is defined \[w'_x(\delta): = \inf_{\{t_i\}}\max_{1\leq i\leq v}\sup_{s,t\in [t_{i-1}, t_i)}|x(s) - x(t)|\]
where $\{t_i\}$ is a partition of $[0,T]$ such that $\displaystyle\max_i\{t_{i}-t_{i-1}\}\leq \delta$ and $x\in D([0,T], \mathbb{R}^+)$ is distributed according to $P_n$.

First, note that since $\mu^{m,n}_t$ is a probability measure, we have
\[|\inprod{f}{\mu^{m,n}_t}|\leq \|f\|_{\infty}\]
for all  $t$, $m$, and $n$. Thus, (i) holds.

For (ii), we have by Markov's inequality:
\begin{align}
\label{eqn_det_tightness}
P_{m,n}(w'(\delta)\geq \epsilon)\leq\tfrac{1}{\epsilon}\mathbb{E}_{m,n}(w'(\delta))
\end{align}
where $w'(\delta):=w'_{\inprod{f}{\mu^{m,n}_t}}(\delta)$. We will use the fact that $\inprod{f}{\mu^{m,n}_t}$ is a pure jump process to bound the right-hand side. The process $\inprod{f}{\mu^{m,n}_t}$ has two types of jumps: nearest-neighbor, and occupied-site jumps. Nearest-neighbor jumps occur at rate
\begin{align*}
\sum_{i} m\mu^{m,n}_t(\tfrac{i}{n})i(1-\tfrac{i}{n})(2+s)\leq \tfrac{mn}{4}(2+s)
\end{align*}
and have magnitude
\begin{align*}
|\inprod{f}{\mu^{m,n}_t}-\inprod{f}{\mu^{m,n}_{t^-}}|& = \left|\left\langle f,\mu^{m,n}_{t^-}+\tfrac{1}{m}\(\delta_{\frac{i\pm 1}{n}}-\delta_{\frac{i}{n}}\)\right\rangle-\inprod{f}{\mu^{m,n}_{t^-}}\right|
\leq \tfrac{1}{mn}\max_{i}|D^-_xf(\tfrac{i}{n})|
\end{align*}
Occupied-site jumps occur at rate
\begin{align*}
\sum_{i,j} m\mu^{m,n}_t(\tfrac{i}{n})\mu^{m,n}_t(\tfrac{j}{n})(1+r\tfrac{j}{n})\leq m(1+r)
\end{align*}
and have magnitude
\begin{align*}
|\inprod{f}{\mu^{m,n}_t}-\inprod{f}{\mu^{m,n}_{t^-}}|& = \left|\left\langle f, \mu^{m,n}_{t^-}+\tfrac{1}{m}\(\delta_{\frac{j}{n}}-\delta_{\frac{i}{n}}\)\right\rangle-\inprod{f}{\mu^{m,n}_{t^-}}\right|\leq \tfrac{2}{m}\|f\|_{\infty}
\end{align*}

Putting this together,
\begin{align*}
\mathbb{E}_{m,n}(w'(\delta))&\leq \mathbb{E}_{m,n}[\text{number of nearest-neighbor jumps in time }\delta] \cdot \tfrac{1}{mn}\max_{i}|D^-_xf(\tfrac{i}{n})|\\
&\qquad + \mathbb{E}_{m,n}[\text{number of occupied-site jumps in time }\delta] \cdot 2\tfrac{1}{m}\|f\|_{\infty}\\
&\leq \tfrac{mn}{4}(2+s)\delta\tfrac{1}{mn}\max_{i}|D^-_xf(\tfrac{i}{n})|+ m(1+r)\delta\tfrac{2}{m}\|f\|_{\infty}\\
&=\left\{\tfrac{2+s}{4}\max_{i}|D^-_xf(\tfrac{i}{n})|+2(1+r)\|f\|_{\infty}\right\}\delta
\end{align*}
Because $f\in C^1([0,1])$, the expression in curly brackets is uniformly bounded by $C_{f}$, a constant that depends on $f$ but not on $m$ nor $n$. Substituting the above into \eqref{eqn_det_tightness} we get that for $\delta<\frac{\epsilon\eta}{C_{f}}$,
\[P_{m,n}(w'(\delta)\geq \epsilon)\leq \eta\]
for all $m$ and $n$. Thus, both conditions for tightness are satisfied and $\inprod{f}{\mu^{m,n}_t}$ is tight.
\end{proof}

\begin{lemma}\label{det_uniqueness}
The integro-partial differential equation \eqref{det_pde} in Theorem 1 has a unique solution.
\end{lemma}

\begin{proof}
Suppose $\mu_t$ satisfies \eqref{det_pde}.
Fix $t\geq 0$ and let $\psi_t(x)$ be a function of time $t$ and space $x$. By the chain rule and the differential equation \eqref{det_pde}, 
\begin{align}
\notag \tfrac{d}{dt}\inprod{\psi_t}{\mu_t}& = \tfrac{d}{dz}\inprod{\psi_z}{\mu_t}\big|_{z=t}+\tfrac{d}{dz}\inprod{\psi_t}{\mu_z}\big|_{z=t}\\
\notag & = \inprod{\fpart{}{t}\psi_t}{\mu_t} -\inprod{sx(1-x)\fpart{\psi_t}{x}}{\mu_t}+ wr\left[\inprod{x\psi_t}{\mu_t}-\inprod{\psi_t}{\mu_t}\inprod{x}{\mu_t}\right]\\
\label{time_dep}\inprod{\psi_t}{\mu_t}&= \inprod{\psi_0}{\mu_0} +\int_0^t \inprod{\fpart{}{z}\psi_z(x) +G\psi_z(x)}{ \mu_z} dz\\
\notag&\qquad\qquad+ wr\int_0^t \inprod{x\psi_z}{\mu_z}-\inprod{\psi_z}{\mu_z}\inprod{x}{\mu_z}dz
\end{align}
where $Gf = -sx(1-x)\fpart{}{x}f$. Let $P_t$ be the semigroup operator associated with  $G$. In fact, using the method of characteristics (or Lemma~\ref{method_char} with $\lambda=0$),
\begin{align}
\label{sgp_op}P_tf = f\(\tfrac{xe^{-st}}{1-x+xe^{-st}}\)
\end{align}

Now, set $\psi_z(x) = P_{t-z}f(x)$ for $0\leq z\leq t$, where $f\in C^1([0,1])$ is some test function. Substituting this into \eqref{time_dep}, we have
\begin{align}
\notag\inprod{P_0f}{\mu_t}=& \inprod{P_t f}{\mu_0} + \int_0^t \inprod{\fpart{}{z}P_{t-z}f(x) +GP_{t-z}f(x)}{ \mu_z} dz\\
\notag & \qquad + \int_0^t wr\left[\inprod{xP_{t-z}f}{\mu_z}- \inprod{P_{t-z}f}{\mu_z}\inprod{x}{\mu_z}\right]dz\\
\label{time_dep_de}\inprod{f}{\mu_t}=& \inprod{P_tf}{\mu_0} + \int_0^t wr\left[\inprod{xP_{t-z}f}{\mu_z}- \inprod{P_{t-z}f}{\mu_z}\inprod{x}{\mu_z}\right]dz
\end{align}
since $\fpart{}{z}P_{t-z}f = -GP_{t-z}f$. Thus, any $\mu_t$ that satisfies \eqref{det_pde} also satisfies \eqref{time_dep_de}. We show that \eqref{time_dep_de} has a unique solution, which in turn implies that \eqref{det_pde} has a unique solution.

Suppose $\mu_t$ and $\nu_t$ both satisfy \eqref{time_dep_de}, with $\mu_0=\nu_0$. Let $t\geq 0$. 
\begin{align}
\notag\|\mu_t-\nu_t\|_{TV}&= \sup_{\|f\|_{\infty}\leq 1}\inprod{f}{\mu_t}-\inprod{f}{\nu_t}\\
&= \sup_{\|f\|_{\infty}\leq 1}\Big\{\int_0^t wr\inprod{xP_{t-z}f}{\mu_z-\nu_z}+ wr\left[\inprod{x}{\mu_z}\inprod{P_{t-z}f}{\mu_z}-\inprod{x}{\nu_z}\inprod{P_{t-z}f}{\nu_z}\right]dz\Big\}\label{tv}
\end{align}
We can bound the first term in the integrand by
\begin{align*}
wr|\inprod{xP_{t-z}f}{\mu_z-\nu_z}|&\leq wr\|\mu_z-\nu_z\|_{TV}
\end{align*}
because $\|xP_{t-z}f\|_{\infty}\leq \|P_{t-z}f\|_{\infty}\leq \|f\|_{\infty}\leq 1$, where the first inequality follows from $x\in[0,1]$ and the second from \eqref{sgp_op}.
For the second term in the integrand of \eqref{tv}, add and subtract $\inprod{x}{\nu_z}\inprod{P_{t-z}f}{\mu_z}$:
\begin{align*}
wr\Big|\inprod{x}{\mu_z}\inprod{P_{t-z}f}{\mu_z}-\inprod{x}{\nu_z}\inprod{P_{t-z}f}{\nu_z}\Big|=& wr\Big|\inprod{x}{\mu_z-\nu_z}\inprod{P_{t-z}f}{\mu_z}+\inprod{x}{\nu_z}\inprod{P_{t-z}f}{\mu_z-\nu_z}\Big|\\
\leq & wr\( \|P_{t-z}f\|_{\infty} \|\mu_z-\nu_z\|_{TV} + \|\mu_z-\nu_z\|_{TV}\)\\
\leq & wr(\|f\|_{\infty}+1)\|\mu_z-\nu_z\|_{TV}
\end{align*}
again, the inequalities follow from $x\in[0,1]$, $\|P_tf\|_{\infty}\leq \|f\|_{\infty}$ and also that $\mu_z$ and $\nu_z$ are probability measures.
Substituting this back into \eqref{tv},
\begin{align*}
\|\mu_t-\nu_t\|_{TV} &\leq \int_0^t 3wr \|\mu_z-\nu_z\|_{TV}dz
\end{align*}
By Gronwall's inequality, $\|\mu_t-\nu_t\|_{TV} = 0$, so we have uniqueness.
\end{proof}

\begin{proof}[Proof of Theorem \ref{det_thm}]
The uniqueness of the limit is given by Lemma~\ref{det_uniqueness} and the tightness of the process by Lemma~\ref{det_tightness}. It remains to show that $\{\inprod{f}{\mu^{m,n}_t}\}_{m,n}$  converges to the solution of \eqref{det_pde}. Recall from Lemma \ref{semi_mart_lem} that
\begin{align*}
 \inprod{f}{\mu^{m,n}_t}-\inprod{f}{\mu^{m,n}_0}& =  A^{m,n}_t(f) + M^{m,n}_t(f)
\end{align*}
Since tightness implies relative compactness (Prohorov's theorem), there exists a subsequence of $\mu^{m,n}_t$ that converges to a limit, call it $\mu_t$. Thus, $\inprod{f}{\mu^{m,n}_t}\to\inprod{f}{\mu_t}$. We also have $\inprod{f}{\mu^{m,n}_0}\to\inprod{f}{\mu_0}$ by assumption. In addition,
\begin{align*}
 A^{m,n}_t(f) & =  \int_0^t\left\{ \sum_i \mu^{m,n}_z(\tfrac{i}{n})\tfrac{i}{n}(1-\tfrac{i}{n})\left[\tfrac{1}{n}D_{xx}f(\tfrac{i}{n}) -sD^-_{x}f(\tfrac{i}{n})\right]\right.\\
&\qquad+ \left.wr\left[\sum_{j}\mu^{m,n}_z(\tfrac{j}{n})\tfrac{j}{n}f(\tfrac{j}{n})-\sum_{i}\mu^{m,n}_z(\tfrac{i}{n})f(\tfrac{i}{n})\sum_j\mu^{m,n}_z(\tfrac{j}{n})\tfrac{j}{n}\right]\right\}~ dz\\
&\rightarrow \int_0^t\left\{\inprod{ -x(1-x)s\tfrac{df}{dx}}{\mu_z}+ wr\left[\inprod{xf(x)}{\mu_z}-\inprod{f(x)}{\mu_z}\inprod{x}{\mu_z}\right]\right\}dz\\
&=:A_t(f)
\end{align*}
The factor of $\frac{1}{m}$ in the quadratic variation \eqref{quad_var} implies that $M^{m,n}_t\to 0$ as $m,n\to\infty$. Therefore, 
\[\inprod{f}{\mu_t}-\inprod{f}{\mu_0} =  A_t(f)\]
or,
\begin{align*}
\tfrac{d}{dt}\inprod{f}{\mu_t}&=\inprod{ -x(1-x)s\tfrac{df}{dx}}{\mu_t}+ wr\left[\inprod{xf(x)}{\mu_t}-\inprod{f(x)}{\mu_t}\inprod{x}{\mu_t}\right]
\end{align*}

\end{proof}

\subsection{Proof of Fleming-Viot limit}\label{stoch_limit}
The elementary proof for tightness in Theorem~\ref{det_thm} does not easily carry over for the case of Theorem~\ref{fv_thm}. We thus use a criterion by Aldous \cite{Aldous:1978} to prove tightness for the martingale part of the stochastic process.



First, consider the semimartingale formulation of $\inprod{f}{\mu^{m,n}_t}$ \eqref{semimart} with the rescaled parameters $s=\frac{\sigma}{n}$ and $\rho = \frac{r}{m}$.
Let $E^{m,n}_t(f):=\int_0^t e^{m,n}_z (f)dz$ and $N^{m,n}_t(f)$ denote the drift and martingale parts of  $\inprod{f}{\nu^{m,n}_t}$, the rescaled process. Then
\begin{align}
\label{rescale_drift}E^{m,n}_t(f) & = \int_0^t\sum_i \nu^{m,n}_z(\tfrac{i}{n})\tfrac{i}{n}(1-\tfrac{i}{n})\left[D_{xx}f(\tfrac{i}{n}) -\sigma D^-_{x}f(\tfrac{i}{n})\right]\\
\notag&\qquad+ w \rho \tfrac{n}{m}\left\{\sum_{j}\nu^{n}_z(\tfrac{j}{n})\tfrac{j}{n}f(\tfrac{j}{n})-\sum_{i}\nu^{n}_z(\tfrac{i}{n})f(\tfrac{i}{n})\sum_j\nu^{n}_z(\tfrac{j}{n})\tfrac{j}{n}\right\} dz
\end{align}
and
\begin{align}
\label{rescale_qv}\langle N^{m,n}(f)\rangle_t=&\int_0^t\left\{\tfrac{n}{m^2}\sum_i \nu^{m,n}_z(\tfrac{i}{n})\tfrac{i}{n}(1-\tfrac{i}{n})\left[\(D^+_xf(\tfrac{i}{n})\)^2+(1+\tfrac{\sigma}{n})\(D^-_xf(\tfrac{i}{n})\)^2\right]\right.\\
\notag& \qquad+\left.w\tfrac{n}{m}\sum_{i,j}\nu^{m,n}_z(\tfrac{i}{n})\nu^{m,n}_t(\tfrac{j}{n})(1+\tfrac{\rho}{m}\tfrac{j}{n})(f(\tfrac{i}{n})-f(\tfrac{j}{n}))^2\right\}dz
\end{align}

\begin{lemma}\label{fv_tightness}
The processes $\inprod{f}{\nu^{m,n}_t}$, as a sequence in $\{(m,n)\}$, is tight for all $f\in C^2([0,1])$.
\end{lemma}
\begin{proof}
Since $\inprod{f}{\nu^{m,n}_t}=E^{m,n}_t(f) + N^{m,n}_t(f)$, it suffices, by the triangle inequality applied to Billingsley's tightness criterion (Theorem~13.2 in \cite{Billingsley:1999}), to show tightness of $E^{m,n}(f)$ and $N^{m,n}(f)$ separately.


For the tightness of the finite variation term $E^{m,n}_t(f)$:
\begin{align*}
|e^{m,n}_t(f)| & \leq \tfrac{1}{4}\sum_i \nu^{m,n}_z(\tfrac{i}{n})\left[\left|D_{xx}f(\tfrac{i}{n})\right| +\sigma \left|D^-_{x}f(\tfrac{i}{n})\right|\right]\\
 &\qquad+ w\rho \tfrac{n}{m}\left\{\sum_{j}\nu^{n}_z(\tfrac{j}{n})\tfrac{j}{n}|f(\tfrac{j}{n})|+\sum_{i}\nu^{n}_z(\tfrac{i}{n})|f(\tfrac{i}{n})|\sum_j\nu^{n}_z(\tfrac{j}{n})\tfrac{j}{n}\right\}
\end{align*}
For a given $\gamma>0$, we can choose $n$ and $m$ sufficiently large such that $\frac{n}{m}\in(\theta-\gamma, \theta+\gamma)$, $|D_{xx}f(\frac{i}{n})|\leq \|f''\|_{\infty}+\gamma$, and $|D^-_xf(\frac{i}{n})|\leq \|f'\|_{\infty}+\gamma$. We thus obtain 
\begin{align*}
|e^{m,n}_t(f)| & \leq \tfrac{1}{4}\left[\|f''\|_{\infty}+\gamma+\sigma(\|f'\|_{\infty}+\gamma)\right]+ 2w\rho(\theta+\gamma) \|f\|_{\infty}
\end{align*}
There are only a finite number of $m$ and $n$ for which this condition is not satisfied. Taking the maximum of the right-hand side of the above equation with the value of $|e^{m,n}_t(f)|$ for such $m$ and $n$, we obtain that for all $m$ and $n$,
\[|e^{m,n}_t(f)|\leq G_{f}\]
and therefore
\[\sup_{t\in [0,T]}|E^{m,n}_t(f)|\leq G_{f}T\]
where $G_{f}$ is a constant that depends on $f$.
Using the same conditions for tightness as in the proof of Theorem~\ref{det_thm}, condition (i) is satisfied because $E^{m,n}_t(f)$ is bounded uniformly in $t$, $m$, and $n$. Condition (ii) is satisfied because $|E^{m,n}_{t+\delta}-E^{m,n}_{t}|\leq \delta G_{f}$ for all $t$, $m$, and $n$ and therefore we can always choose $\delta$ to be sufficiently small so that $|E_{t+\delta}^{m,n}-E_{t}^{m,n}|\leq \epsilon$ for some prescribed $\epsilon$.

We will show tightness for the martingale part $\langle N_t^{m,n}(f)\rangle_t$ using Aldous' tightness condition (we use the result as stated in \cite{Etheridge:2000}). First, note that by equation~\eqref{rescale_qv},
\begin{align*}
\langle N^{m,n}_t(f)\rangle_t \leq J_{f} t
\end{align*}
for $f\in C^2([0,1])$, where $J_{f}$ is a constant that depends on $f$. Thus for fixed $t$,
\begin{align*}
P_{m,n}(|N^{m,n}_t(f)|>a)&\leq \tfrac{1}{a}\mathbb{E}_{m,n}|N^{m,n}_t(f)|\\
&\leq \tfrac{1}{a} \(\mathbb{E}_{m,n}(N^{m,n}_t(f))^2\)^{1/2}=\tfrac{1}{a}\(\mathbb{E}_{m,n}\langle N^{m,n}_t(f)\rangle_t\)^{1/2}\leq\tfrac{\sqrt{J_{f}t}}{a}
\end{align*}
Given $\epsilon>0$, choose $a>\frac{\sqrt{J_{f}t}}{\epsilon}$ and we have that $N^{m,n}_t(f)$ is tight for each $t$.
Next, let $\tau$ be a stopping time, bounded by $T$, and let $\epsilon>0$. For $\kappa>0$,
\begin{align*}
P_{m,n}(|N^{m,n}_{\tau+\kappa}(f)-N^{m,n}_{\tau}(f)|\geq \epsilon)&\leq \tfrac{1}{\epsilon}\mathbb{E}_{m,n}|N^{m,n}_{\tau+\kappa}(f)-N^{m,n}_{\tau}(f)|
\end{align*}
Now (suppressing subscripts on expected value for clarity),
\begin{align*}
\mathbb{E}|N^{m,n}_{\tau+\kappa}(f)-N^{m,n}_{\tau}(f)|&\leq \left[\mathbb{E}(N^{m,n}_{\tau+\kappa}(f)-N^{m,n}_{\tau}(f))^2\right]^{1/2}\\
&=\left[\mathbb{E}(N^{m,n}_{\tau+\kappa}(f)^2-N^{m,n}_{\tau}(f)^2+2N^{m,n}_{\tau}(f)(N^{m,n}_{\tau}(f)-N^{m,n}_{\tau+\kappa}(f))\right]^{1/2}\\
&=\left[\mathbb{E}(\langle N^{m,n}(f)\rangle_{\tau+\kappa}-\langle N^{m,n}(f)\rangle_{\tau})\right]^{1/2}\\
&\leq \sqrt{J_{f}\kappa}
\end{align*}
Hence,
\begin{align*}
P_{m,n}(|N^{m,n}_{\tau+\kappa}(f)-N^{m,n}_{\tau}(f)|\geq \epsilon)&\leq \tfrac{1}{\epsilon}\sqrt{J_{f}\kappa}
\end{align*}
By taking $\kappa<\frac{\epsilon^4}{\sqrt{J_{f}}}$, we satisfy the conditions of Aldous' stopping criterion.

\end{proof}

\begin{lemma}\label{fv_uniqueness}
The martingale problem \eqref{fv_drift} and \eqref{fv_quad} has a unique solution.
\end{lemma}
\begin{proof}
The martingale problem with $V(t,\nu, x)=0$ corresponds to a neutral Fleming-Viot with linear mutation operator. Its uniqueness has previously been established (see for example  \cite{Dawson:2010}). To show uniqueness for nontrivial $V$, we use a Girsanov-type transform by Dawson \cite{Dawson:2010}. It suffices to check that 
\begin{align}
\sup_{t,\mu,x} |V(t, \mu, x)|\leq V_0 \text{ (a constant)}
\end{align}
In our case, $V(t, \mu, x)=x$ and since $x\in[0,1]$, the condition is satisfied and the martingale problem has a unique solution.
\end{proof}

\begin{proof}[Proof of Theorem \ref{fv_thm}]
The uniqueness of the limit is given by Lemma~\ref{fv_uniqueness} and the tightness of the process by Lemma~\ref{fv_tightness}. To see that the limit is the martingale problem stated in Theorem \ref{fv_thm}, note that for a fixed $t$,
\begin{align*}
E^{m,n}_t(f) &\longrightarrow \int_0^t \int_0^1x(1-x)[\spart{}{x}f(x) - \sigma\fpart{}{s}f(x)]\nu_z(dx)\\
&\qquad\qquad w\rho\theta \left\{\int_0^1xf(x)\nu_z(dx)-\int_0^1f(x)\nu_z(dx)\int_0^1x\nu_z(dx)\right\}dz
\end{align*}
as $n,m\rightarrow\infty$ and
\begin{align*}
\langle N^{m,n}(f)\rangle_t \longrightarrow \int_0^tw\theta \int_0^1\int_0^1(f(x)-f(y))^2\nu_z(dx)\nu_z(dy)dz
\end{align*}
Finally, notice that 
\begin{multline*}
\int_0^1\int_0^1(f(x)-f(y))^2\nu_z(dx)\nu_z(dy)\\
= 2\int_0^1\int_0^1f(x)^2\nu_z(dx) \nu_z(dy)- 2\int_0^1\int_0^1f(x)f(y)\nu_z(dx)\nu_z(dy)\\
= 2\int_0^1\int_0^1 f(x)f(y)\nu_z(dx)[\delta_x(dy)-\nu_z(dy)]
\end{multline*}
and
\begin{multline*}
\int_0^1xf(x)\nu_z(dx)-\int_0^1f(x)\nu_z(dx)\int_0^1x\nu_z(dx)
= \int_0^1\int_0^1f(x)y\nu_z(dx)[\delta_x(dy)-\nu_z(dx))]
\end{multline*}
satisfying the form of the martingale problem in the theorem.
\end{proof}

\newpage
\bibliographystyle{plain}

	\bibliography{weak_convergence_bib.bib}

\begin{thebibliography}{10}

\bibitem{Aldous:1978}
David Aldous.
\newblock Stopping times and tightness.
\newblock {\em The Annals of Probability}, 6(2):pp. 335--340, 1978.

\bibitem{Billingsley:1999}
Patrick Billingsley.
\newblock {\em Convergence of Probability Measures}.
\newblock Wiley-Interscience, 2 edition, 1999.

\bibitem{Champagnat:2006}
Nicolas Champagnat, R\'{e}gis Ferri\`{e}re, and Sylvie M\'{e}l\'{e}ard.
\newblock Unifying evolutionary dynamics: From individual stochastic processes
  to macroscopic models.
\newblock {\em Theor Pop Biol}, 69(3):297--321, 2006.

\bibitem{Dawson:2010}
Donald~A Dawson.
\newblock {Introductory Leture on Stochastic Population Systems, Technical
  Report Series of the Laboratory for Research Statistics and Probability}.
\newblock Technical Report 451, Carleton University - University of Ottawa,
  2010.

\bibitem{Dawson:1991}
Donald~A Dawson and K~J Hochberg.
\newblock A multilevel branching model.
\newblock {\em Advances in Applied Probability}, 23:701--715, 1991.

\bibitem{Durrett:2008}
Richard Durrett.
\newblock {\em Probability Models for {DNA} Sequence Evolution}.
\newblock Springer, 2nd edition, 2008.

\bibitem{Etheridge:2000}
Alison Etheridge.
\newblock {\em An Introduction to Superprocesses}.
\newblock American Mathematical Soc., 2000.

\bibitem{Ethier:1993}
S.~N. Ethier and Thomas~G. Kurtz.
\newblock Fleming-viot processes in population genetics.
\newblock 31(2):42, 1993.

\bibitem{Fleming:1979}
Wendell~H. Fleming and Michel Viot.
\newblock Some measure-valued markov processes in population genetics theory.
\newblock {\em Indiana University Mathematics Journal}, 28(5):817--843, 1979.

\bibitem{Fournier:2004}
Nicolas Fournier and Sylvie M\'el\'eard.
\newblock A microscopic probabilistic description of a locally regulated
  population and macroscopic approximations.
\newblock {\em Ann. Appl. Probab.}, 14(4):1880--1919, 2004.

\bibitem{Joffe:1986}
A~Joffe and M~Metivier.
\newblock Weak convergence of sequences of semimartingales with applications to
  multitype branching processes.
\newblock {\em Advances in Applied Probability}, 18:20--65, 1986.

\bibitem{Kallenberg:1997}
Olav Kallenberg.
\newblock {\em Foundations of Modern Probability}.
\newblock Springer, 1st edition, 1997.

\bibitem{Luo:2014}
Shishi Luo.
\newblock A unifying framework reveals key properties of multilevel selection.
\newblock {\em J Theor Biol}, 341:41--52, 2014.

\bibitem{Pinchover:2005}
Yehuda Pinchover and Jacob Rubinstein.
\newblock {\em An Introduction to Partial Differential Equations}.
\newblock Cambridge University Press, 2005.

\bibitem{Protter:2004}
Philip~E. Protter.
\newblock {\em Stochastic Integration and Differential Equations}.
\newblock Springer, 2nd edition, 2004.

\end{thebibliography}

\end{document}